\documentclass[a4paper]{amsart}
\usepackage{amsmath,amssymb,textcomp,graphicx}
\usepackage{mathrsfs}
\usepackage{times}
\usepackage[T1]{fontenc}
\usepackage[english]{babel}
\usepackage{rawfonts}
\usepackage{color}
\usepackage[latin1]{inputenc}
\renewcommand{\d}{\textnormal{d}}

\newcommand{\ii}{\infty }

\newcommand{\R}{\mathbb{R}}

\newcommand{\N}{\mathbb{N}}

\newtheorem{rema}{Remark}[section]
\newtheorem{theorem}{Theorem}[section]
\newtheorem{lemme}{Lemma}[section]

\title{Shape Minimization of Dendritic Attenuation}
\author{Antoine Henrot\and
Yannick Privat}
\address{Institut \'Elie Cartan de Nancy, UMR 7502
Nancy-Universit\'e - CNRS - INRIA, B.P. 239 , Vand\oe
uvre-l\`es-Nancy Cedex France}
\thanks{Corresponding author: Antoine Henrot, Tel
+33-383684560, Fax +33-383684534, E-mail address
Antoine.Henrot@iecn.u-nancy.fr}
\date{\today}

\begin{document}

\begin{abstract}
What is the optimal shape of a dendrite? Of course, optimality refers
to some particular criterion. In this paper, we look at the case of
a dendrite sealed at one end and connected at the other end to a soma.
The electrical potential in the fiber follows the classical cable
equations as established by W. Rall. We are interested in the shape
of the dendrite which minimizes either the attenuation in time of the
potential or the attenuation in space. In both cases, we prove that
the cylindrical shape is optimal.
\end{abstract}

\maketitle
\section{Introduction}
\subsection{Motivation}
Is Nature always looking for optimum for living organisms? In
particular, are the organs designed to optimize some criterion?
Complete answers to these questions are likely never to be
discovered. Nevertheless, assuming that Nature proceeds in the most
efficient way, can lead to a better understanding of the modeling of
an organ and the underlying physical or chemical phenomena. This is
this idea of {\it inverse modeling} that we had in mind when we
began this work. Roughly speaking, it can be described by the
following steps:
\begin{enumerate}
  \item Let us consider a given organ of a living body.
  \item Write a mathematical model which describes the behavior of
  this organ.
  \item Imagine a (numerical) criterion that Nature would like to
  optimize for this organ.
  \item Determine the optimal shape for this criterion and this
  model.
  \item Compare with the real shape(s).
\end{enumerate}
If the optimal and the real shapes coincide, we can guess that our
model and our criterion are relevant. If they do not, we must admit
that either our criterion or our model (or the initial guess that
Nature looks for optimum) is probably wrong. We believe that it will
often be the choice of the criterion which is not correct. A
possible reason is the complexity of Nature. This complexity
indicates that, in general, there is not a unique criterion to
optimize but several ones (which could also be antagonists). The
mathematical study (point 4 in the above procedure) becomes then
much more difficult since one needs to use tools of multi-criteria
optimization.

In this paper, we want to follow the above procedure in the case
of a dendrite. We consider a fiber which is sealed at its right end
and connected to a soma at its left end. We use the classical
cable equation to describe the electrical potential along the
fiber. What are the criterions that we can consider here? Of
course, we want the dendrite to propagate the best as possible the
electrical signal. In other terms, the attenuation of the signal
must be as small as possible. We are going to consider the two
kinds of possible attenuation: {\bf attenuation in time or in
space} and we are looking for the shape of a dendrite which minimizes
this attenuation. In both cases, the optimal shape that we get is
a cylinder. Since it is very close to the real shape, at first
sight, we can conclude that Nature is in accordance with
mathematics for this problem and solves a shape optimization
problem! For a general reference in mathematical modeling in 
Neuroscience, we refer to the book of A. Scott, \cite{scott}. 
For  a more exhaustive view and 
an introduction to the beauty of shapes (and, in
particular, optimal shapes) in Nature, we refer to the classical
books of S. Hildebrandt and A. Tromba, \cite{HT} and A. Bejan, \cite{bej}.
\subsection{The mathematical model}
Let us consider a fiber with a cylindrical symmetry, of length
$\ell$ and radius $a(x)$ at point $x$. We denote by $v(x,t)$ the
difference from rest of the membrane potential at point $x$ and time
$t$. The equation satisfied by $v(x,t)$ is similar to the classical
cable equation as established by W. Rall during the sixties, cf
\cite{rall1}, \cite{rall2}, \cite{cox}. See also, \cite{st-sp} as the
best motivation of the current study. We consider here the case of
a fiber which is sealed at its right end and connected to a soma
with surface area $A_s$ at its left end. Let us denote by $R_a$ the
axial resistance (k$\Omega$cm), $C_m$ is the membrane capacitance
($\mu$F/cm$^2$), $G_m$ the fiber membrane conductance and $G_s$ the
soma membrane conductance (mS/cm$^2$). We assume that the fiber is
initially at rest and that it receives a transient current stimulus
$i_0$ at the left. The parabolic equation satisfied by $v$ is then
(see \cite{cox}):
\begin{equation}\label{1}
\left\{\begin{array}{lr}\vspace{2mm} \frac{1}{2R_a}\frac{\partial
}{\partial x}\left(a^2\frac{\partial v}{\partial x}\right)
=a\sqrt{1+a'^2}\left(C_m\frac{\partial v}{\partial t}+G_mv\right)&
x\in (0,\ell),\,t>0
\\\vspace{2mm}
\frac{\pi a^2(0)}{R_a}\frac{\partial v}{\partial x}(0,t)=A_s
\left(C_m\frac{\partial v}{\partial t}(0,t)+G_sv(0,t)\right)-i_0(t)& t>0 \\
\frac{\partial v}{\partial x}(\ell , t )=0 & t>0\\
v(x,0)=0 & x\in (0,\ell).
\end{array}
\right.
\end{equation}
It is convenient to represent $v$, solution of (\ref{1}), in terms
of eigenfunctions as did S. Cox and J. Raol in \cite{cox}:
\begin{equation}
\label{DecompSoma} v(x,t)=\sum_{n=0}^{+\infty}\psi _n(t)\phi _n(x)
\ \quad  x \in (0,\ell ), \  t> 0
\end{equation}
where $\phi_n$ is the $n$-th eigenfunction associated to the
eigenvalue $\mu _n$ :
\begin{equation}\label{2}
\left\{\begin{array}{ll}\vspace{2mm} -(a^2\phi _n')'=\mu _n\,
a\sqrt{1+a'^2}\phi _n & x\in (0,\ell )\\\vspace{2mm}
\frac{2\pi}{A_s} \, a^2(0)\phi _n'(0)+(\mu_n+\gamma )\phi _n(0)=0 & \\
\phi _n'(\ell )=0 & \\
\end{array}
\right.
\end{equation}
where $\gamma :=2R_a(G_m-G_s)$ is assumed to be non negative. We
choose to normalize the eigenfunctions by
\begin{equation}\label{norma}
\Arrowvert \phi_n \Arrowvert_a^2
 := A\phi_n^2(0)+\int_0^\ell
a(x)\sqrt{1+a'^2(x)}\phi_n^2(x)\d x=1.
\end{equation}
where $A=\frac{A_s}{2\pi}$. Of course the eigen-pair
$(\mu_n,\phi_n)$ strongly depends on the taper $a(x)$ so we
will often denote it by $(\mu_n(a),\phi_n(a))$. The eigenvalue
problem (\ref{2}) is not classical since the eigenvalue appears in
the boundary condition, see section \ref{TimeAttenuation} for more
precisions.
\par Inserting the decomposition (\ref{DecompSoma}) in the equation
(\ref{1}) gives an ordinary differential equation satisfied by
$\psi_n(t)$. After resolution, we get:
\begin{equation}\label{soluv}
v(x,t)=\frac{1}{2\pi C_m}\sum_{n=1}^{+\infty}\phi _n(0) \phi _n(x
)\, i_0* e^{-\lambda_n t}\, ,
\end{equation}
where $\lambda _n :=\frac{\mu_n+2R_aG_m}{2R_aC_m}>0$ (see Lemma
\ref{lemmasign}) and $*$ denotes the convolution product of distributions.

\subsection{The optimization problems}
We need now to give a precise statement to the optimization problems
presented in the introduction. In that purpose, we have to choose
the functions we want to optimize and the class of functions $a(x)$
in competition. Let us begin with this last point. Since the fiber
must not collapse, it is natural to assume a lower bound for the
functions $a(x)$, so we fix a positive constant $a_0$ and we impose:
\begin{equation}\label{lower}
    \forall x\in [0,\ell],\quad a(x)\geq a_0>0\,.
\end{equation}
Now the minimal regularity needed for $a$ is clearly, according to
system (\ref{1}) or (\ref{2}) that the derivative $a'$ exists (at
least almost everywhere) and is bounded, so we choose to work in the
class of Lipschitz continuous functions which is often denoted by
$W^{1,\infty}(0,\ell )$. At last, we also need to put a constraint
on the "cost" for Nature to build a fiber. It seems reasonable to
consider that this cost is proportional to the surface area of the
fiber. This surface area is clearly given by
\begin{equation}\label{surface}
    \displaystyle \mbox{Surface area}=2\pi \int_0^{\ell}a(x)
    \sqrt{1+a'^2(x)} \d x
\end{equation}
so we can assume a bound, say $S$, on this surface area. To
summarize, we consider the class of functions $a(x)$ defined by:
\begin{equation}\label{class}
\mathcal{A}_{a_0 , S}:=\left\{a\in W^{1,\infty}(0,\ell ), \ a(x)
\geq a_0 \textrm{ and }\int_0^{\ell}a(x)\sqrt{1+a'^2(x)}\d x \leq
S\right\}.
\end{equation}
Of course, we need to assume $S> a_0\ell$ in order that the
class $\mathcal{A}_{a_0 , S}$ be non trivial.

\medskip
As explained in the Introduction, an "ideal" dendrite should conduct in
the best possible way the electrical information he is supposed to
transmit. In other terms, the attenuation of the signal must be
minimized. Since the potential $v$ depends on the space and the time
variable, we can consider both criterions.
\subsubsection{Attenuation in space}
Let us introduce the transfer function $T$ defined by :
\begin{equation}\label{transfer}
    T(a):=\frac{\int_0^{+\infty}v(0,t)\d t}{\int_0^{+\infty}v(\ell ,t)\d
t}.
\end{equation}
$T(a)$ corresponds to the ratio of the mean values in time of the
potential $v$ taken at points $x=0$ and $x=\ell$. This ratio is
always greater than one, see Remark \ref{rem24} and is a good
indicator of the attenuation of the signal between the two
extremities of the dendrite. So, it is realistic to look for a taper
$a(x)$ which yields a ratio as close to one as possible:
\begin{equation}\label{optimspace}
\textsl{Find $a$ in the class $\mathcal{A}_{a_0, S}$ which minimizes
$T(a)$.}
\end{equation}
\subsubsection{Attenuation in time}
According to expansion (\ref{soluv}), the potential $v(x,t)$ goes to
0 when $t\to +\ii$. More precisely, its asymptotic behavior is
described by
$$v(x,t)\simeq \frac{1}{2\pi
C_m} \phi _1(0) \phi _1(x ) \, i_0* e^{-\lambda_1 t}$$
where $\lambda _1 :=\frac{\mu_1(a)+2R_aG_m}{2R_aC_m}>0$ and
$\mu_1(a)$ is the first eigenvalue of (\ref{2}). Therefore, as it is
classical in such parabolic problems, it seems natural to look for a
function $a(x)$ which minimizes the exponential rate of decay:
\begin{equation}\label{optimtime}
\textsl{Find $a\in \mathcal{A}_{a_0, S}$ which minimizes $\mu _1(a)$
(the first eigenvalue of (\ref{2})).}
\end{equation}
The idea of minimizing eigenvalues of such Sturm-Liouville operators
is a long story and goes back at least to M. Krein in \cite{krein},
see also \cite{henrot} for a review on such problems.
\subsubsection{The main result}
We state in the following Theorem the main results of this paper
\begin{theorem}
\label{maintheo}
\begin{enumerate}
\item The unique minimizer of the eigenvalue $\mu_1 (a)$ in the class
$\mathcal{A}_{a_0,S}$ is the constant function $a\equiv a_0$.
\item The unique minimizer of the criterion $T (a)$ in the class
$\mathcal{A}_{a_0,S}$ is the constant function $a\equiv a_0$.
\end{enumerate}
\end{theorem}
In other terms, for both criterions, the optimal shape of a dendrite
sealed at one end and connected to a soma at the other end is the
cylindrical one!

\medskip In his thesis and in a foregoing paper, see \cite{privat1},
\cite{privat2}, the second author studies the case of a dendrite
sealed at both ends. From a mathematical point of view, it changes
the boundary conditions in (\ref{1}) and (\ref{2}) which become
homogeneous Neumann boundary conditions at both extremities. The
result he obtains is the same for the case of attenuation in space,
but it is different for the attenuation in time. Actually, there is
no existence of a minimizer for $\mu_1(a)$ (as usual in the Neumann
case $\mu_0(a)=0$ and $\mu_1(a)$ denotes the first non-zero
eigenvalue). Moreover, he is able to exhibit minimizing sequences
which would produce very strange dendrites!

\subsection{Notation}
The set of notation used in this paper is summarized in this
section.
\begin{center}
\begin{tabular}{ll}\vspace{3mm}
$W^{1,\ii}(0,\ell)$ & the set of Lipschitz continuous functions
defined on the interval $[0,\ell]$.\\\vspace{3mm} $\mathcal{A}_{a_0,
S}$ & \begin{tabular}{ll} the class of functions
defined by\\
$\left\{a\in W^{1,\infty}(0,\ell), \ a(x) \geq a_0,\;
\int_0^{\ell}a(x)\sqrt{1+a'^2(x)}\d x \leq S\right\}$.
                          \end{tabular}
\\\vspace{3mm}
$\Arrowvert . \Arrowvert _\infty$ & \begin{tabular}{ll} norm defined
on the space of bounded functions $L^\infty(0, \ell)$ by\\
$\Arrowvert f \Arrowvert_\infty =\sup _{t\in [0, \ell]}|f(t)|$.
\end{tabular}
\\\vspace{3mm} $<.,.>_a$ & \begin{tabular}{ll} inner product defined for two continuous
functions $f$ and $g$ by:\\
$<f,g>_a := Af(0)g(0)+\int_0^\ell
a(x)\sqrt{1+a'^2(x)}f(x)g(x)\d x.$
\end{tabular}
\\\vspace{3mm}
$\Arrowvert . \Arrowvert_a$ & norm induced by
$<.,.>_a$.\\\vspace{3mm} $\mathcal{E}_a$ & completion of the space
of continuous functions $\mathcal{C}([0,\ell ])$ for the norm
$\Arrowvert . \Arrowvert _a$.\\\vspace{3mm} $L^2(0,\ell)$ & the
space of (classes of) functions which are square-integrable on
$(0,\ell)$.\\\vspace{3mm} $H^1(0,\ell)$ & \begin{tabular}{ll} the
Sobolev space of functions in $L^2(0,\ell)$ whose derivative (in the
sense\\ of distributions) lies in $L^2(0,\ell)$.
\end{tabular}
\\\vspace{3mm} $
\left< \frac{\d J}{\d \nu }(\nu _0),h\right>$ &
\begin{tabular}{ll} G\^ateaux-derivative of a function $J$ at point
$\nu _0$ in direction $h$ defined by:\\ $ \left< \frac{\d J}{\d \nu
}(\nu _0),h\right>\stackrel{\textnormal{def}}{=} \lim_{t\searrow
0}\frac{J(\nu _0+t h)-J (\nu _0)}{t}. $\end{tabular}
\end{tabular}
\end{center}

\section{Minimization of the first eigenvalue}
\label{TimeAttenuation} The eigenvalue problem (\ref{2}) is not
completely classical due to the presence of the eigenvalue $\mu_n$
in the first boundary condition. As explained in \cite{cox}, see
also the works of J. Walter \cite{walter} and J. Ercolano-M.
Schechter \cite{ES}, a good way to handle with such case consists in
introducing the following inner product:
$$<f,g>_a = Af(0)g(0)+
\int_0^\ell a(x)\sqrt{1+a'^2(x)}f(x)g(x)\d x,$$ its associated norm
$\|.\|_a$ and the Hilbert space $\mathcal{E}_a$ defined as the
completion of the space of continuous functions $\mathcal{C}([0,\ell
])$ for this norm. It is easy to see that $\mathcal{E}_a$ is a space
satisfying $H^1(0,\ell)\subset \mathcal{E}_a \subset L^2(0,\ell)$
(both inclusions are strict). Moreover, the map $\phi\mapsto
\phi(0)$ defines a linear continuous form on $\mathcal{E}_a$. It is
now classical spectral theory which allows to prove existence of a
sequence of eigenvalues $\mu_n$ and eigenfunctions $\phi_n$
orthogonal for the inner product $<.,.>_a$.

\medskip
Let us now make an elementary observation on the sign of $\mu_1(a)$:
\begin{lemme}\label{lemmasign}
Let $a$ be in the class $\mathcal{A}_{a_0, S}$, then the first
eigenvalue $\mu_1(a)$ of (\ref{2}) satisfies
$$-2R_aG_m< -\gamma < \mu_1(a) < 0\,.$$
\end{lemme}
\begin{proof}
Let $a\in \mathcal{A}_{a_0,S}$ and $v\in H^1(0, \ell)$, non
identically zero. We denote by $R[a;v]$, the Rayleigh quotient:
$$
R[a;v]:=\frac{\displaystyle \int_0^{\ell}a^2(x)u'^2(x) \d x-A\gamma
u^2(0)}{\displaystyle \int_0^{\ell}a(x)\sqrt{1+a'^2(x)}u^2(x) \d
x+Au^2(0)}.
$$
The classical Poincar\'e-Courant-Hilbert principle writes:
\begin{equation}\label{PCH}
    \mu _1(a)=\inf_{v\in H^1(0,\ell)} R[a;v]
\end{equation}
Now, taking $v \equiv 1$ in the above formula yields:
$$
\mu_1 (a) \leq -\frac{A\gamma}{\displaystyle
\int_0^{\ell}a(x)\sqrt{1+a'^2(x)} \d x+A} < 0 \,.
$$
The lower bound is easy by observing that $\mu_1(a)+\gamma=
R[a;\phi_1(a)]+\gamma >0$. 
\end{proof}
\begin{rema} \rm
Using the min-max formulae for the second eigenfunction, it is also
possible to prove that the second eigenvalue satisfies $\mu_2(a)>0$.
The proof consists in studying the problem of calculus of variations
$\min  \frac{1}{\displaystyle v^2(0)}.\int_0^{\ell}a^2(x)v'^2(x)\d
x$ on the class $W_a:=\{v\in H^1(0,\ell ) : v(0)\neq 0\textrm{ and
}<v,\phi_1(a)>_{a}=0\}$. \end{rema}

We can now prove the first part of Theorem \ref{maintheo}. The
eigenfunction $\phi_1(a)$ associated to $\mu_1(a)$ realizes the
minimum of the Rayleigh quotient. Hence, we have:
$$
\mu _1(a)=\frac{\displaystyle \int_0^{\ell}a^2(x){\phi_1(a)}'^2(x)
\d x-A\gamma {\phi_1(a)}^2(0)}{\displaystyle
\int_0^{\ell}a(x)\sqrt{1+a'^2(x)}{\phi_1(a)}^2(x) \d
x+A{\phi_1(a)}^2(0)}.
$$
\par By Lemma \ref{lemmasign}, the numerator of this quotient is negative.
Moreover, we have:
\begin{equation}
\label{rayleigh1} \int_0^{\ell}a^2(x){\phi_1(a)}'^2(x) \d x-A\gamma
{\phi_1(a)}^2(0)\geq
a_0^2\int_0^{\ell}{\phi_1(a)}'^2(x) \d x-A\gamma {\phi_1(a)}^2(0)
\end{equation}
and \begin{equation}\label{rayleigh2}
\frac{1}{\int_0^{\ell}a\sqrt{1+a'^2}\phi_1(a)(x)\d
x+A{\phi_1(a)}^2(0)} \leq \frac{1}{\int_0^{\ell}a_0\phi_1(a)(x)\d
x+A{\phi_1(a)}^2(0)}
\end{equation}
(with a strict inequality in (\ref{rayleigh1}), (\ref{rayleigh2}) if
$a$ is not constant). Finally, writing that $\mu
_1(a_0)=\displaystyle \inf_{v\in H^1(0,\ell)} R[a_0;v]$, we deduce
from (\ref{rayleigh1}), (\ref{rayleigh2}) and (\ref{PCH}) that $\mu
_1(a_0)< \mu _1(a)$ as soon as $a\not= a_0$.

\section{Minimization of $T(a)$}
\subsection{Introduction}
This last section is devoted to the proof of the second claim of
Theorem \ref{maintheo}, i.e. that $a=a_0$ minimizes the criterion
$T(a)$. We recall that the criterion $T(a)$ describes the
attenuation in space and that it is defined by
$$
    T(a):=\frac{\int_0^{+\infty}v(0,t)\d t}{\int_0^{+\infty}v(\ell ,t)\d
t}.
$$
The proof here is much more complicated than for $\mu_1(a)$. Let us now
outline the different steps of the proof.
\begin{description}
  \item[1st step] Using the Laplace Transform $\widehat{v}(x,p)$ of $v(x,t)$,
  we rewrite the criterion $T(a)$ as the quotient
  $\widehat{v}(0,0)/\widehat{v}(\ell,0)$.
  \item[2nd step] We use the change of variable defined by $y=\int_0^x \frac{\d
  t}{a^2(t)}$ to transform the equation satisfied by $\widehat{v}$ into
  a simpler differential equation. This allows us to consider a new
  unknown $\rho(y):=a^3(x)\sqrt{1+a'(x)^2}$ instead of $a$ and a new
  criterion $T_1(\rho)$. The function $\rho$ must lie in the set
  defined by:
$$\mathcal{R}_{a_0,S}:=\left\{\rho \in L^{\infty}(0, \ell _1) :
\ a_0^3\leq \rho (y)\textrm{ and }\int_0^{\ell _1}\rho (y)\d y\leq S
\right\}.$$
  \item[3rd step] We solve the new optimization problem $\min
  T_1(\rho)$ first on the subclass of functions $\rho\in
  \mathcal{R}_{a_0,S}$ which satisfy $\rho\leq M$ for some positive constant
  $M$. We prove that the
  minimizer has to be a \textsl{bang-bang} function. It means that it can only
  takes the values $a_0^3$ and $M$. Then, studying
  carefully the optimality conditions, we prove that the only
  minimizer is $\rho\equiv a_0^3$.
  \item[4th step] We conclude.
\end{description}
\subsection{Use of the Laplace Transform}
The parabolic equation is not completely standard in the sense that
it contains a dynamical boundary condition at $x=0$. This kind of
problem has been studied by different people, see e.g.
\cite{Escher}, \cite{BVR} and the references therein. It can be
proved that the solution $v(x,t)$ belongs to $L^2(0,T,H^1(0,\ell))$.
Moreover, using the eigenvalue expansion (\ref{DecompSoma}), one can
see that, in the case of an impulsion $i_0=\delta$ (a Dirac measure
at $t=0$), the integrals $\int_0^{+\infty}|v(0,t)|\d t$ and
$\int_0^{+\infty}|v(\ell,t)|\d t$ are well defined. Let us introduce
the Laplace Transform $\widehat{v}(x,p)$ of the solution defined by
$$\widehat{v}(x,p):=\int_0^{+\infty}e^{-pt}v(x,t)\d t\,.$$
Thanks to the convergence of the integrals, the criterion $T(a)$ can
be rewritten as
\begin{equation}\label{3.1}
    T(a):=\frac{\int_0^{+\infty}v(0,t)\d t}{\int_0^{+\infty}v(\ell ,t)\d
t}=\frac{\lim_{p\to 0}\widehat{v}(0,p)}{\lim_{p\to 0}
\widehat{v}(\ell,p)}\,.
\end{equation}
Now, transforming equation (\ref{1}), we see that the Laplace
Transform $\widehat{v}$ is the solution of the following o.d.e.:
\begin{equation}\label{edpLaplace}
\left\{\begin{array}{ll}\vspace{2mm}
\frac{1}{2R_a}\frac{\partial}{\partial x}\left(
a^2\frac{\partial\widehat{v}}{\partial
x}\right)=a\sqrt{1+a'^2}\left(C_m p\widehat{v}+G_m\widehat{v}\right)
& (x,p)\in (0,\ell )\times (0;+\infty)\\\vspace{2mm} \frac{\pi
a^2(0)}{R_a}\frac{\partial \widehat{v}}{\partial x}(0,p)=A_s
\left[C_mp\widehat{v}(0,p)+G_s\widehat{v}(0,p)\right]-1 & p\in (0,+\infty)\\\vspace{2mm}
\frac{\partial \widehat{v}}{\partial x}(\ell , p )=0 & p\in (0,+\infty) \\
\widehat{v}(x,0)=0 & x\in (0, \ell).
\end{array}
\right.
\end{equation}
\subsection{A change of variable}
We are now going to use the following change of variable classical
in ordinary differential equations, see e.g. \cite{coxlip}:
$$
y=\int_0^x \frac{\d t}{a^2(t)}.
$$
The interval $(0,\ell)$ becomes $(0,\ell_1)$ where
$$
    \ell _1:=\displaystyle \int_0^{\ell}\frac{\d t}{a^2(t)},
$$
the function $\widehat{v}$ is transformed into the function
$$
w(y,p):=\widehat{v}(x,p).
$$
and we consider a new unknown $\rho$ defined by :
\begin{equation}\label{defrho}
    \rho (y):=a^3(x)\sqrt{1+{a'}(x)^2}\,.
\end{equation}
Since $a$ belongs to the class $\mathcal{A}_{a_0,S}$ defined in
(\ref{class}), the new function $\rho$ belongs to:
\begin{equation}\label{class2}
    \mathcal{R}_{a_0,S}:=\left\{\rho \in L^{\infty}(0, \ell _1) :
\ a_0^3\leq \rho (y)\textrm{ and }\int_0^{\ell _1}\rho (y)\d y\leq S
\right\}.
\end{equation}
Then, equation (\ref{edpLaplace}) becomes:
\begin{equation}\label{edpCV}
\left\{\begin{array}{ll}\vspace{2mm} \frac{1}{2R_a}\frac{\partial
^2w}{\partial y^2}=\rho \left(C_mp+G_m\right)w & (y,p)\in (0,\ell
_1)\times (0;+\infty)\\ \vspace{2mm} \frac{\pi }{R_a}\frac{\partial
w}{\partial y}(0,p)=A_s
\left[C_mp+G_s\right]w(0,p)-1 & p\in (0,+\infty)\\
\frac{\partial w}{\partial y}(\ell _1, p )=0 & p\in (0,+\infty)\,.
\end{array}
\right.
\end{equation}
We let $p$ going to 0 in the equation (\ref{edpCV}) (see Appendix
\ref{annexA} for a mathematical justification) to get a function
$w_0(y):=w(y,0)$ which satisfies
\begin{equation}\label{edpCVFinal}
\left\{\begin{array}{ll}\vspace{2mm} \frac{1}{2R_a}\frac{\d
^2w_0}{\d y^2}=\rho G_m w_0 & y\in (0,\ell _1)\\ \vspace{2mm}
\frac{\pi }{R_a}\frac{\d w_0}{\d y}(0)=A_s G_s w_0(0)-1 & \\
\frac{\d w_0}{\d y}(\ell _1)=0 & .
\end{array}
\right.
\end{equation}
Moreover, from (\ref{3.1}) the criterion $T(a)$ becomes
\begin{equation}\label{T1}
    T(a)=T_1(\rho)=\frac{w_0(0)}{w_0(\ell _1)}
\end{equation}
The problems $\min \{T(a), \, a\in \mathcal{A}_{a_0,S}\}$ and 
$\min \{T_1(\rho), \, \rho\in \mathcal{R}_{a_0,S}\}$ are not completely
equivalent since $a\in\mathcal{A}_{a_0,S} \mapsto \rho\in\mathcal{R}_{a_0,S}$
is not a one-to-one correspondance. Nevertheless, it is clear 
that $\mathcal{R}_{a_0,S}$ contains the whole image of $\mathcal{A}_{a_0,S}$
by this map. So, if we find a minimizer of $T_1$ in $\mathcal{R}_{a_0,S}$
which belongs to the image of $\mathcal{A}_{a_0,S}$ (this will be the case),
we will solve our problem.
\begin{rema}\label{rem24} \rm
Let us have a look to equation (\ref{edpCVFinal}). It is not
possible that $w'_0(0)\geq 0$ (otherwise $w_0(0)$ would be positive,
according to the first boundary condition, and then $w_0$ would
remain positive and convex which would contradict the second
boundary condition. In the same way, $w_0(0)$ cannot be negative, otherwise
$w_0$ would remain negative and concave and this is impossible with
the second boundary condition.So finally, one can see
that $w'_0(0)<0$, $w_0(0)>0$ and $w_0$ remains positive on the whole
interval $(0,\ell_1)$. At last, since $w'_0$ is increasing and
$w'_0(\ell_1)=0$, we see that $w_0$ is decreasing (and positive),
therefore $w_0(\ell_1)<w_0(0)$ which proves that $T(a)=T_1(\rho)<1$.
\end{rema}
\subsection{Study of a new optimization problem}
Using the different transformations introduced in the previous
subsections, we see that we must now solve the optimization problem:
$\min T_1(\rho)$ with $\rho$ in the class $\mathcal{R}_{a_0,S}$. Let
us begin by solving this optimization problem in the subclass
$$
\mathcal{R}_{a_0,S}^M:=\left\{\rho \in L^{\infty}(0, \ell _1) : \
a_0^3\leq \rho (y)\leq M\textrm{ and }\int_0^{\ell _1}\rho (y)\d
y\leq S \right\}
$$
where $M$ is a fixed positive constant ($M>a_0^3$). We will let
$M\to +\infty$ later.
\begin{theorem}\label{theobang}
The problem $\min T_1(\rho)$, with $\rho\in\mathcal{R}_{a_0,S}^M$,
has a solution $\rho^*$. Moreover, every solution is a {\it
bang-bang} function, i.e. a function which satisfies $\rho^*=a_0^3$
or $M$ almost everywhere.
\end{theorem}
\begin{proof}
The fact that the optimizer is a bang-bang function often occurs in
such control problems. Existence of a minimizer $\rho^*$ is easy,
due to continuity of the criterion $T_1(\rho)$ for the weak-*
convergence. Then, we write and analyze the optimality conditions
thanks to the introduction of two adjoint problems. We are able to
prove that the set $\{a_0^3<\rho^*<M\}$ has zero measure. Let us now
give the details.\\
Since $\mathcal{R}_{a_0,S}^M$ is a bounded subset of $L^{\infty}(0,
\ell _1)$ it is compact for the weak-star convergence . So, to prove
existence of a minimizer, we just need to prove that the criterion
$T_1$ is continuous for the weak-star convergence. Let $(\rho
_n)_{n\in \N}$ be a sequence in $\mathcal{R}_{a_0,S}^M$ such that
$\rho _n \stackrel{*}{\rightharpoonup}\rho$ and let us denote by
$w_n$ and $w$ the associated solutions of (\ref{edpCVFin}). From the
variational formulation of this problem
\begin{equation}\label{fvwn}
\frac{1}{2R_a}\int_0^{\ell _1}w_n'(y)z'(y)\d y+G_m\int_0^{\ell_1}
\rho_n (y)w_n(y)z(y)\d y+2AG_s w_n(0)z(0)= \frac{z(0)}{2\pi}
\end{equation}
for every $z\in H^1(0,\ell_1)$, we first see (taking $z=w_n$ in
(\ref{fvwn}) that the sequence $w_n$ is bounded in $H^1(0,\ell_1)$.
So, it converges (up to a subsequence) weakly in $H^1(0,\ell_1)$ and
strongly in $L^2(0,\ell_1)$ to a function $w_\infty$. Now, these
convergence are sufficient to pass to the limit in (\ref{fvwn}), so
we have proved that $w_\infty=w$ and all the sequence converges
since $w$ is the only accumulation point. Existence of a minimizer
$\rho^*$ in the class $\mathcal{R}_{a_0,S}^M$ follows.

We want now to write the optimality conditions. Let  $h\in
W^{1,\infty}(0, \ell _1)$ be an admissible perturbation of the
optimum. We will now denote by $\dot{w_0}$ the quantity:
$$
\dot{w_0}:=\left<\frac{\d w_0}{\d \rho}(\rho^* ),h\right>.
$$
Classical variational analysis shows that $\dot{w_0}$ is the
solution of the ordinary differential equation:
\begin{equation}\label{wdot}
\left\{\begin{array}{ll}\vspace{2mm} \frac{1}{2R_a}\frac{\d
^2\dot{w_0}}{\d y^2}=G_m \left(\rho \dot{w_0}+hw_0\right) & \quad
y\in (0,\ell _1)\\ \vspace{2mm}
\frac{\pi }{R_a}\frac{\d \dot{w_0}}{\d y}(0)=A_s G_s \dot{w_0}(0) & \\
\frac{\d \dot{w_0}}{\d y}(\ell _1)=0. &
\end{array}
\right.
\end{equation}
Differentiating the criterion $T_1$ in the direction $h$ gives:
\begin{equation}\label{exT1}
\left<\frac{\d T_1}{\d \rho},h\right>= \frac{\dot{w_0}(0)w_0(\ell
_1)-w_0(0)\dot{w_0}(\ell _1)}{w_0^2(\ell _1)}.
\end{equation}
Let us now introduce the two following adjoint problems. We
consider the function $q_1$ defined as the solution of the ordinary
differential equation:
\begin{equation}\label{q1}
\left\{\begin{array}{ll}\vspace{2mm} \frac{1}{2R_a}\frac{\d
^2q_1}{\d y^2}=G_m \rho \left(q_1(y)-y\right) & y\in (0,\ell _1)\\
\vspace{2mm}
\frac{\pi }{R_a}\frac{\d q_1}{\d y}(0)=A_s G_s q_1(0) & \\
\frac{\d q_1}{\d y}(\ell _1)=0 & 
\end{array}
\right.
\end{equation}
and the function $q_2$ solution of :
\begin{equation}\label{q2}
\left\{\begin{array}{ll}\vspace{2mm} \frac{1}{2R_a}\frac{\d
^2q_2}{\d y^2}=G_m \rho \left(q_2(y)-1\right) & y\in (0,\ell _1)\\
\vspace{2mm}
\frac{\pi }{R_a}\frac{\d q_2}{\d y}(0)=A_s G_s q_2(0) & \\
\frac{\d q_2}{\d y}(\ell _1)=0 & 
\end{array}
\right.
\end{equation}
(existence and uniqueness of $q_1$, $q_2$ follows from Lax-Milgram
Theorem). Multiplying equation (\ref{wdot}) by $q_2-1$, equation
(\ref{q2}) by $\dot{w_0}$ and integrating both by parts yields:
\begin{equation}\label{wdot0}
    \dot{w_0}(0)=\frac{G_m}{AG_s}\,\int_0^{\ell_1} h(y) w_0(y)
    (q_2(y)-1)\,dy\;.
\end{equation}
In the same way,  multiplying equation (\ref{wdot}) by $q_1-y$,
equation (\ref{q1}) by $\dot{w_0}$ and integrating both by parts
yields:
\begin{equation}\label{wdotl}
    \dot{w_0}(\ell_1)=\dot{w_0}(0)+2R_aG_m\int_0^{\ell_1} h(y) w_0(y)
    (q_1(y)-y)\,dy\;.
\end{equation}
Therefore, (\ref{exT1}) together with (\ref{wdot0}) and
(\ref{wdotl}) gives
\begin{equation}
\left<\frac{\d T_1}{\d \rho},h\right> =
\frac{2R_aG_m}{w_0(\ell_1)^2}\, \int _0^{\ell _1}h(y)w_0(y)
\left(\frac{w_0(\ell
_1)-w_0(0)}{\tilde{A}}\left(q_2(y)-1\right)-w_0(0)
\left(q_1(y)-y\right)\right)\d y
\end{equation}
where $\tilde{A}=2AG_sR_a$. Let us denote by $f$, the function of
one variable defined by:
$$
f:\begin{array}[t]{rcl}
[0, \ell _1] & \longrightarrow & \R \\
y & \longmapsto &
\frac{2R_aG_m}{w_0(\ell_1)^2}\,\left\lbrack\frac{w_0(\ell
_1)-w_0(0)}{\tilde{A}}
\left(q_2(y)-1\right)-w_0(0)\left(q_1(y)-y\right)\right\rbrack.
\end{array}
$$
So, finally:
$$
\left<\frac{\d T_1}{\d \rho},h\right>= \int _0^{\ell
_1}h(y)w_0(y)f(y)\d y.
$$
We want now to prove that the optimum $\rho^*$ is a bang-bang
function. It is a classical approach, see e.g.
\cite{henrot-maillot}. For that purpose, let us introduce the
following sets:
\begin{itemize}
\item $\mathcal{I}_0(\rho^* ):=\left\{y\in (0, \ell _1) : \rho^* (y)=a_0\right\}$ ;
\item $\mathcal{I}_M(\rho^* ):=\left\{y\in (0, \ell _1) : \rho^* (y)=M\right\}$ ;
\item $\mathcal{I}_\star(\rho^* ):=(\left\{y\in (0, \ell _1) : a_0<\rho^*
(y)<M\right\}$.
\end{itemize}
We write  $\mathcal{I}_\star(\rho^* )= \displaystyle \bigcup
_{k=1}^{+\infty}\left\{y\in (0, \ell _1) : a_0+\frac{1}{k}< \rho^*
(y)<M-\frac{1}{k}\right\}=\bigcup _{k=1}^{+\infty}\mathcal{I}_{\star
, k}(\rho^*)$. We want to prove that $\mathcal{I}_{\star , k}(\rho^* )$ has
zero measure, for all integer $k\neq 0$. We argue by contradiction.
Let us suppose that $|\mathcal{I}_{\star , k}(\rho^* )|\neq 0$. Let
$y_0\in \mathcal{I}_{\star,k}(\rho^* )$. We denote by
$(G_{k,n})_{n\geq 0}$, the sequence of subsets of
$\mathcal{I}_{\star , k}$ :
$$
G_{k,n}:=B\left(y_0,\frac{1}{n}\right)\cap \mathcal{I}_{\star ,
k}(\rho^*)\subset \mathcal{I}_{\star , k}(\rho^*).
$$
Let us notice that $\displaystyle
\bigcap_{n=0}^{+\infty}G_{k,n}=\{y_0\}$, and let us choose $h=\chi
_{G_{k,n}}$. Then, for $t$ small enough, perturbations $ \rho^* +th$
et $\rho^* -th$ are admissible. Then:
$$
\lim_{t\searrow 0}\frac{T_1 (\rho^*
+th)-T_1(\rho^*)}{t}=\int_{0}^\ell h(y)w_0(y)f(y)\d y\geq
0\Longleftrightarrow \int_{G_{k,n}} w_0(y)f(y)\d y\geq 0.
$$
In the same way:
$$
\lim_{t\searrow 0}\frac{T_1 (\rho^*
-th)-T_1(\rho^*)}{t}=-\int_{0}^\ell h(y)w_0(y)f(y)\d y\geq
0\Longleftrightarrow \int_{G_{k,n}} w_0(y)f(y)\d y\leq 0.
$$
We can deduce that $\displaystyle \int_{G_{k,n}} w_0(y)f(y)\d y= 0$.
We divide by $|G_{k,n}|$ and we make $n$ tending to $+\infty$.
The Lebesgue density theorem shows that $w(y_0)f(y_0)=0$,
a.e. for $y_0\in \mathcal{I}_{\star ,k}(\rho^* )$. This is clearly a
contradiction, since $w_0$ and $f$ are respectively solutions of the
differential equations $\frac{1}{2R_a}\frac{\d ^2w_0}{\d y^2}=\rho^*
G_m w_0$ and $\frac{1}{2R_a}\frac{\d ^2f}{\d y^2}=\rho^* G_m f$.
This proves that $|\mathcal{I}_{\star , k}(\rho^* )|=0$ and then
$\mathcal{I}_\star(\rho^* )$ has also zero measure, what implies
that $\rho^*$ equals $a_0^3$ or $M$ almost everywhere.
\end{proof}

\medskip\noindent
Now, we prove that among every bang-bang function, this is the
constant function $a_0^3$ which yields the minimum of $T_1$.
\begin{lemme}\label{lemtechn}
The optimum of $T_1$ in the class $\mathcal{R}_{a_0,S}^M$ is the
constant function $\rho(y)=a_0^3$.
\end{lemme}
\begin{proof}
Using notation of the proof of Theorem \ref{theobang}, the
optimality conditions writes
\begin{itemize}
\item On the set $\mathcal{I}_0(\rho^* )$, we have $f(y)\geq 0$ and $h(y)\geq 0$ ;
\item On the set $\mathcal{I}_M(\rho^* )$, we have $f(y)\leq 0$ and $h(y) \leq 0$.
\end{itemize}
According to the differential equation (\ref{q2}) and maximum
principle, the function $q_2-1$ is negative on $[0, \ell _1]$. Let
us write $f$ like below:
$$
f(y)=\frac{2R_aG_m}{w_0(\ell_1)^2}\,\left(q_2(y)-1\right).\left(\frac{w_0(\ell
_1)-w_0(0)}{\tilde{A}}-w_0(0)g(y)\right),
$$
where $\displaystyle g(y):=\frac{q_1(y)-y}{q_2(y)-1}$. The
function $g$ is two times derivable on $(0, \ell _1)$, and, using
equations (\ref{q1}), (\ref{q2}) we have:
\begin{equation}
\forall y\in (0, \ell _1), \ \frac{\d ^2g}{\d
y^2}(y)=-2\frac{\frac{\d }{\d y}[q_2-1](y)}{q_2(y)-1}.\frac{\d g}{\d
y}(y).
\end{equation}
Therefore, on every interval where $\frac{\d g}{\d y}$ keeps his
sign, we have:
$$
\frac{\d g}{\d y}(y)=\frac{C}{[q_2(y)-1]^2}, \textrm{ with }C\in \R
.
$$
Now, $\frac{\d g}{\d y}$ is continuous on $[0, \ell _1]$, then the
only possibility is:
\begin{equation}
\forall y \in (0, \ell _1), \ \frac{\d g}{\d y}(y)=-\frac{q_2(\ell
_1) -1}{[q_2(y)-1]^2}>0.
\end{equation}
It follows that $g$ is an increasing function on $[0, \ell
_1]$. Since $q_2-1$ is negative and $w_0(0)>0$, we see that: $f(y)\geq
0\Longleftrightarrow g(y)\geq \frac{w(\ell
_1)-w(0)}{\tilde{A}w(0)}$. Then, according to the optimality
conditions the function $\rho^*$, local optimum for the criterion
$T_1$, is necessarily as follows:
\begin{equation}\label{loc1}
    \rho^*(y)=\left\lbrace
    \begin{array}{cc}
      M & \mbox{if } y<\xi_1 \\
      a_0^3 & \mbox{if } y>\xi_1 \\
    \end{array}\right.
\end{equation}
with a transition point $\xi_1$ which is possibly 0 or $\ell_1$.
\begin{figure}[h]\label{bangbang}
\begin{center}
\includegraphics[width=10cm,height=8cm]{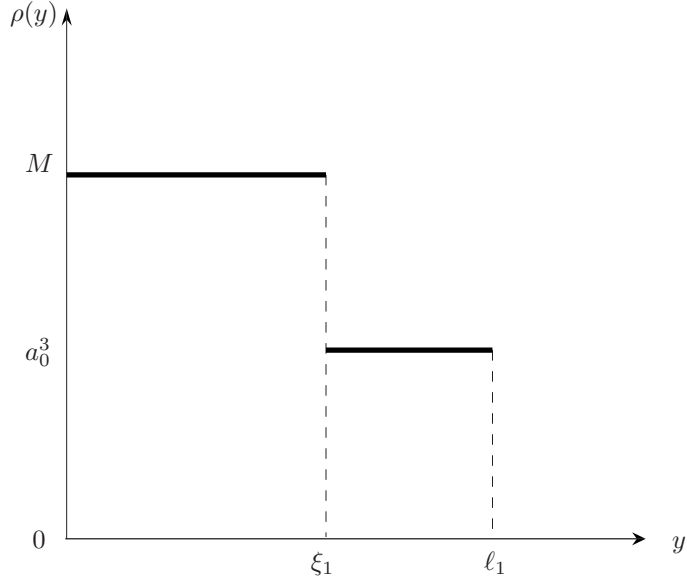}
\caption{Possible profile of the optimum}
\begin{picture}(0,0)
\put(35,20){$y$}
\put(-51,70){$M$}
\put(-51,45){$a_0^3$}
\put(-53,90){$\rho (y)$}
\put(-50,20){$0$}
\put(-13,17){$\xi_1$}
\put(10,17){$\ell_1$}
\end{picture}
\end{center}
\end{figure}
Now, it remains to look for the better $\rho^*$ among all functions
defined by (\ref{loc1}). The only unknown is finally the transition
point $\xi_1$. We can use the explicit expression of the solution
$w_0$ of equation (\ref{edpCVFin}) for such a simple $\rho^*$.
We write $w_0$
\begin{equation}
w_0 (y)=\left\{\begin{array}{ll}
\alpha _1^M\cosh (\omega _My)+\alpha _2^M\sinh (\omega _My) & \textrm{ on }[0,\xi_1]\\
\alpha _3^M\cosh (\omega _0y)+\alpha _4^M\sinh (\omega _0y) &
\textrm{ on }[\xi_1 , \ell _1]
\end{array}
\right.
\end{equation}
where $\alpha _1^M$, $\alpha _2^M$, $\alpha _3^M$ and $\alpha _4^M$
are four constants and $\omega _0:=\sqrt{2R_aG_ma_0^3}$, $\omega
_M:=\sqrt{2R_aG_mM}$. Thanks to boundary conditions, we get
\begin{equation}
w_0 (y)=\left\{\begin{array}{ll}
\alpha _1^M\cosh (\omega _My)+
\frac{A_sG_s\alpha _1^M-1}{\frac{\pi}{R_a}\omega _M}\sinh (\omega _My) &
\textrm{ on }[0,\xi_1]\\
\alpha _3^M\left(\cosh (\omega _0y)-\tanh (\omega _0\ell _1 )\sinh
(\omega _0y)\right) & \textrm{ on }[\xi_1 , \ell _1]
\end{array}
\right.
\end{equation}
Therefore, the criterion $T_1$ is given by
$$
T_1(\rho^*
)=\frac{w_0(0)}{w_0(\ell_1)}\,=\frac{\alpha_1^M}{\alpha_3^M}\cosh
(\omega_0 \ell _1).
$$
Finally, using the fact that $w_0$ is $C^1$, we get, thanks to
continuity of $w_0$ and $\frac{d w_0}{d y}$ at $y=\xi_1$
\begin{eqnarray}
\frac{\alpha_1^M}{\alpha_3^M} & = & \cosh (\omega _0\xi_1)\cosh
(\omega _M\xi_1)-
\frac{\omega _0}{\omega _M}\sinh (\omega _M\xi_1)\sinh (\omega _0\xi_1) \nonumber \\
 & & -\tanh (\omega _0\ell _1)\left[\sinh (\omega _0\xi_1)\cosh (\omega _M\xi_1)-
 \frac{\omega _0}{\omega _M}\sinh (\omega _M\xi_1)\cosh (\omega _0\xi_1)\right]
\end{eqnarray}
and then
\begin{eqnarray}
T_1(\rho^* ) & = & \cosh (\omega _0\ell _1)\left[\cosh (\omega
_0\xi_1)\cosh
(\omega _M\xi_1)-\frac{\omega _0}{\omega _M}\sinh (\omega _M\xi_1)\sinh (\omega _0\xi_1)\right] \nonumber \\
 & & -\sinh (\omega _0\ell _1)\left[\sinh (\omega _0\xi_1)\cosh (\omega _M\xi_1)-
 \frac{\omega _0}{\omega _M}\sinh (\omega _M\xi_1)\cosh (\omega _0\xi_1)\right]
\end{eqnarray}
Computing the derivative of the criterion with respect to the
variable $\xi_1$, we have
\begin{eqnarray}
\frac{\d T_1(\rho^* )}{\d \xi_1} & = & \left(\frac{\omega
_M^2-\omega _0^2}{\omega _M}\right) \sinh (\omega _M \xi_1)\cosh
[\omega _0(\ell _1-\xi_1)] \geq 0.
\end{eqnarray}
It follows that $\xi_1$ has to be equal to 0, that means that the
constant function $a_0^3$ minimize the criterion $T_1$ on
$\mathcal{R}_{a_0,S}^M$.
\end{proof}

\subsection{Conclusion}
The proof of Theorem \ref{maintheo} follows now easily. Since
$a_0^3$ is the unique minimizer of the criterion  $T_1$ in the class
$\mathcal{R}_{a_0,S}^M$ and
$$
\mathcal{R}_{a_0,S}=\bigcup_{M>a_0^3}\mathcal{R}_{a_0,S}^M.
$$
we get that $a_0^3$ is the (unique) minimizer of $T_1$ in the class
$\mathcal{R}_{a_0,S}$. Moreover, since
$a_0^3=a_0^3\sqrt{1+\left(\frac{\d a_0}{\d y}\right)^2}$, it is
clear that $T_1(a_0^3)=T(a_0)$ and then, $a_0$ minimizes $T$ on the
class $\mathcal{A}_{a_0,S}$.
\appendix
\section{Limit of $w_p$ when $p\to 0$}\label{annexA}
Let us now denote by $w_p$ the function $w(.,p)$ for a given
positive $p$. The equation (\ref{edpCV}) can be rewritten:
\begin{equation}\label{edpCVbis}
\left\{\begin{array}{ll}\vspace{2mm} \frac{1}{2R_a}\frac{\d
^2w_p}{\d y^2}=\rho \left(C_mp+G_m\right)w_p & y\in (0,\ell _1)\\
\vspace{1mm}
\frac{\pi }{R_a}\frac{\d w_p}{\d y}(0)=A_s \left[C_mp+G_s\right]w_p(0)-1 & \\
\frac{\d w_p}{\d y}(\ell _1)=0 \,.& 
\end{array}
\right.
\end{equation}
We recall that we want to prove that $w_p$ has a limit $w_0$ when
$p\to 0$ and that $w_0$ is the solution of (\ref{edpCVFinal}). We can
suppose that $p\in (0,1]$. Let us write the variational formulation
of (\ref{edpCVbis}):
$$\forall z\in H^1(0, \ell _1), \ a_p(w_p,z)=L(z),
$$
with
$$
a_p(u,z)  =  \frac{1}{2R_a}\int_0^{\ell _1}u'(y)z'(y)+(C_mp+G_m)
\rho (y)u(y)z(y)\d y+A(C_mp+G_s)u(0)z(0)$$ and $$L(z) =
\frac{z(0)}{2\pi}.
$$
Let us make $z=w_p$ in the above formulation. Since $a_p$ is
uniformly (with respect to $p$) coercive:
\begin{equation}\label{unicor}
    \min \left(\frac{1}{2R_a},a_0^3 G_m\right)\Arrowvert w_p\Arrowvert
^2_{H^1(0, \ell _1)}\leq a_p(w_p,w_p)=L(w_p)
\end{equation}
and $L$ is a linear continuous form, we get from
(\ref{unicor}) that the sequence $w_p$ is bounded in
$H^1(0,\ell_1)$. Therefore, there exists $w^*\in H^1(0,\ell _1)$
such that $\displaystyle w_p \underset{p\to 0}{\rightharpoonup} w^*$
in $H^1(0, \ell _1)$ and $w_p\xrightarrow[p\to 0]{} w^*$ in $L^2(0,
\ell _1)$ (up to a subsequence). It follows that
\begin{itemize}
  \item $w_p(0)\xrightarrow[p\to 0]{} w^*(0)$,
  \item $\displaystyle \int
_0^{\ell _1}\rho (y)w_p(y)z(y) \d y\xrightarrow[p\to 0]{}\int
_0^{\ell _1}\rho (y)w^*(y)z(y) \d y$,
\item $
\displaystyle \int _0^{\ell _1} w_p'(y)z(y) \d y \xrightarrow[p \to
0]{} \int _0^{\ell _1} w^{*}{'} (y)z'(y)\d y. $
\end{itemize}
Therefore $w^*$ is the solution of the ordinary differential
equation:
\begin{equation}\label{edpCVFin}
\left\{\begin{array}{ll}\vspace{2mm} \frac{1}{2R_a}\frac{\d
^2w_0}{\d y^2}=\rho G_m w_0 & y\in (0,\ell _1)\\ \vspace{2mm}
\frac{\pi }{R_a}\frac{\d w_0}{\d y}(0)=A_s G_s w_0(0)-1 & \\
\frac{\d w_0}{\d y}(\ell _1)=0 \,.& 
\end{array}
\right.
\end{equation}
which means that $w^*=w_0$, the desired result.





\end{document}